\documentclass[12pt]{amsart}
\usepackage[all]{xy}
\usepackage{graphicx}
\usepackage{amsmath,amsxtra,amssymb,latexsym, amscd, amsthm, amscd}
\usepackage{indentfirst}
\usepackage[mathscr]{eucal}

\newtheorem{thm}{Theorem}[section]
\newtheorem{cor}[thm]{Corollary}
\newtheorem{lem}[thm]{Lemma}

\theoremstyle{definition}

\newtheorem{exm}[thm]{Example}
\newtheorem{rem}[thm]{Remark}

\DeclareMathOperator{\lk}{lk}

\DeclareMathOperator{\pd}{pd}
\DeclareMathOperator{\ordmatch}{order-match}
\DeclareMathOperator{\match}{match}

\DeclareMathOperator{\reg}{reg}

\DeclareMathOperator{\Min}{Min} 
 
\DeclareMathOperator{\smp}{\mathcal{SP}}
\DeclareMathOperator{\N}{\mathbb {N}}
\DeclareMathOperator{\Z}{\mathbb {Z}}
\DeclareMathOperator{\R}{\mathbb {R}}
\DeclareMathOperator{\simp}{Simp}

\def\mi {\mathfrak m}
\def\ni {\mathfrak n}
\def\P {\mathcal P}

\def\h {\widetilde{H}}

\def\alb {\boldsymbol {\alpha}}
\def\btb {\boldsymbol {\beta}}

\def\gmb {\boldsymbol {\gamma}}

\def\F{{\mathcal F}}
\def\P{{\mathcal P}}
\def\C{{\mathcal C}}
\def\H{{\mathcal H}}
\def\G{{\mathcal G}}

\def\h{{\widetilde H}}

\def\v {\mathbf v}
\def\x {\mathbf x}
\def\y {\mathbf y}

\begin{document}

\title[Regularity of symbolic powers of square-free monomial ideals]{Regularity of symbolic powers of square-free monomial ideals}

\author[T.T. Hien]{Truong Thi Hien}
\address{Faculty of Natural Sciences, Hong Duc University, No. 565 Quang Trung, Dong Ve, Thanh Hoa, Vietnam}
\email{hientruong86@gmail.com}

\author[T.N. Trung]{ Tran Nam Trung}
\address{Institute of Mathematics, VAST, 18 Hoang Quoc Viet, Hanoi, Viet Nam, and Institute of Mathematics and 
TIMAS, Thang Long University, Ha Noi, Vietnam.}
\email{tntrung@math.ac.vn}

\date{}
\dedicatory{}
\commby{}
\begin{abstract}

We study the regularity of symbolic powers of square-free monomial ideals. We prove that if $I = I_\Delta$ is the Stanley-Reisner ideal of a simplicial complex $\Delta$, then $\reg(I^{(n)}) \leqslant \delta(n-1) +b$ for all $n\geqslant 1$, where $\delta = \lim\limits_{n\to\infty} \reg(I^{(n)})/n$, and $b = \max\{\reg(I_\Gamma) \mid  \Gamma \text{ is a subcomplex of } \Delta \text{ with }  \F(\Gamma) \subseteq \F(\Delta)\}$. This bound is sharp for any $n$. When $I = I(G)$ is the edge ideal of a simple graph $G$, we obtain a general linear upper bound
$\reg(I^{(n)}) \leqslant 2n + \ordmatch(G)-1$, where  $\ordmatch(G)$ is the ordered matching number of $G$.
\end{abstract}
\maketitle
\section*{Introduction}

Throughout the paper, let $K$ be a field and $R = K[x_1,\ldots,x_r]$ the polynomial ring of $r$ variables $x_1,\ldots,x_r$ with $r \geqslant 1$. Let $I$ be a homogeneous ideal of $R$. Then the $n$-th symbolic power of $I$ is defined by
$$I^{(n)} = \bigcap_{\mathfrak{p}\in \Min(I)} I^n R_{\mathfrak{p}} \cap R,$$ 
where $\Min(I)$ is as usual the set of minimal associated prime ideals of $I$.

Cutkosky, Herzog, Trung \cite{CHT}, and independently Kodiyalam \cite{K}, proved that the function $\reg(I^n)$ is a linear function in $n$ for $n\gg 0$. The similar result for symbolic powers is not true even when $I$ is a square-free monomial ideal (see e.g. \cite[Theorem 5.15]{DHHT}) except for the case $\dim(R/I)\leqslant 2$ (see \cite{HT2}).

If $I$ is a square-free monomial ideal, Hoa and the second author (see \cite[Theorem 4.9]{HT1}) proved that the limit
\begin{equation}\label{EQ-Limit}
\delta(I) = \lim\limits_{n\to\infty} \frac{\reg (I^{(n)})}{n},
\end{equation}
does exist, in fact the limit exists for arbitrary monomial ideals (see \cite{DHHT}). Moreover, $\reg (I^{(n)}) < \delta(I)n +\dim(R/I)+1$ for all $n\geqslant 1$. This bound is obvious not sharp for every $n$ (see Corollary \ref{main-cor}). There have been many recent results which establish sharp bounds for $\reg (I^{(n)})$ in the case $I$ is the edge ideal of a simple graph (see e.g. \cite{BN, F, GHOS, JNS}).

The aim of this paper is to find sharp bounds for $\reg (I^{(n)})$, for a square-free monomial ideal $I$, in terms of combinatorial data from its associated simplicial complexes and hypergraphs.

\medskip

For a simplicial complex $\Delta$ on the set $V=\{1,\ldots,r\}$, the Stanley-Reisner ideal of $\Delta$ is defined by
$$I_{\Delta} = \left(\prod_{i\in \tau} x_i \mid \tau \subseteq V \text{ and } \tau \notin \Delta
\right) \subseteq R.$$
Let us denote by $\F(\Delta)$ the set of all facets of $\Delta$.

\medskip

The first main result of the paper is the following theorem.

\medskip

\noindent{\bf Theorem $\ref{main-thm}$.} {\it Let $\Delta$ be a simplicial complex. Then,
$$\reg(I_\Delta^{(n)}) \leqslant \delta(I_\Delta)(n -1)+ b, \ \text{ for all } n\geqslant 1,$$
where 
$b = \max\{\reg(I_\Gamma) \mid  \Gamma \text{ is a subcomplex of } \Delta \text{ with }  \F(\Gamma) \subseteq \F(\Delta)\}$.}
\medskip

This bound is sharp for every $n$ (see Example \ref{exam-matroid}). It is worth mentioning that the number $\delta(I_\Delta)$, which is determined by Equation (\ref{EQ-Limit}), may be not an integer and even bigger than $\reg(I_\Delta)$ (see \cite[Lemma 5.14 and Theorem 5.15]{DHHT}).
\medskip

\medskip

For a simple hypergraph $\H=(V,E)$ with vertex set $V=\{1,\ldots,r\}$, the edge ideal of $\H$ is defined by
$$I(\H) = \left(\prod_{i\in e} x_i \mid e\in E\right)\subseteq R.$$

Let $\H^*$ be the simple hypergraph corresponding to the Alexander duality $I(\H)^*$ of $I(\H)$. Let $\epsilon(\H^*)$ be the minimum number of cardinality of edgewise dominant sets of $\H^*$, this concept was introduced by Dao and Schweig \cite{DS}. 
\medskip

Then second main result of the paper is the following theorem.

\medskip

\noindent{\bf Theorem $\ref{main-hyper}$.} {\it Let $\H$ be a simple hypergraph. Then, 
$$\reg(I(\H)^{(n)}) \leqslant \delta(I(\H))(n -1)+ |V(\H)| - \epsilon(\H^*), \ \text{ for all } n\geqslant 1.$$}

A hypergraph is a graph if every edge has exactly two vertices. For a graph $G$,  a linear lower bound for $\reg(I(G)^{(n)})$ is given in \cite{GHOS}:
$$\reg(I(G)^{(n)}) \geqslant 2n + \nu(G)-1,$$
where $\nu(G)$ is the induced matching number of $G$.  Note that this lower bound is also valid for ordinary powers (see \cite[Theorem 4.5]{BHT}).

On the upper bounds, Fakhari (see \cite[Conjecture 1.3]{F}) conjectured that
$$\reg(I(G)^{(n)}) \leqslant 2n + \reg(I(G))-2,$$
This conjecture, it may be the best bound up to now of our knowledge.

By using Theorem \ref{main-thm}, we obtain a general linear upper bound  for $\reg(I(G)^{(n)})$ in terms of the ordered matching number of $G$, although it is weaker than the one in this conjecture, it provides   us a sharp bound. Note that this result also settles  the question (2) of Fakhari in \cite{Fc}.

\medskip

\noindent{\bf Theorem $\ref{EdgeBoundLinear}$.} { \it Let $G$ be a graph. Then,
$$\reg(I(G)^{(n)}) \leqslant 2n + \ordmatch(G)-1, \text{ for all } n\geqslant 1,$$
where $\ordmatch(G)$ is the ordered matching number of $G$.}

\medskip

Let us explain the idea to prove Theorems \ref{main-thm} and \ref{main-hyper} as follows. Let $i\geqslant 0$ such that $\reg(R/I^{(n)}) = a_i(R/I^{(n)})+i$. 

The first key point is to prove that $a_i(R/I^{(n)}) \leqslant \delta(I)(n -1)$. Assume that $\alb = (\alpha_1,\ldots,\alpha_r) \in \Z^r$ such that
$$H_\mi^i(R/I^{(n)})_{\alb} \ne 0, \text{ and } a_i(R/I^{(n)}) = |\alb|,$$
where $\mi = (x_1,\ldots,x_r)$ and $|\alb| = \alpha_1+\cdots+\alpha_r$.  We reduce to the case $\alb \in \N^r$. In order to bound $|\alb|$, we use Takayama's formula (see Lemma \ref{TA}) to compute $H_\mi^i(R/I^{(n)})_{\alb}$, which allows us to search for $\alb$ in a polytope in $\R^r$, so that we can get the desired bound of $|\alb|$ via theory of convex polytopes (see Theorem \ref{a-Invariant}). 

The second key point is to bound the index $i$ by using the regularity of a Stanley-Reisner ideal in terms of the vanishing of reduced homology of simplicial complexes which derived from Hochster's formula about the Hilbert series of  the local cohomology module of Stanley-Reisner ideals (see Lemma \ref{Hochster-Reg}).

\medskip

Our paper is structured as follows. In the next section, we collect notations and terminology used in the paper, and recall a few auxiliary results. In Section $2$, we prove  Theorems \ref{main-thm} and \ref{main-hyper}. In the last section, we prove Theorem \ref{EdgeBoundLinear}.

\section{Preliminaries}  We shall follow standard notations and terminology from usual texts in the research area (cf. \cite{E,HH2, MS}). For simplicity, we denote the set $\{1,\ldots,r\}$ by $[r]$.

\subsection{Regularity and projective dimension}
Through out this paper, let $K$ be a field, and let $R = K[x_1,\ldots,x_r]$ be a standard graded polynomial ring of $r$ variables over $K$. The object of our work is the Castelnuovo-Mumford regularity of graded modules and ideals over $R$. This invariant can be defined via either the minimal free resolutions or the local cohomology modules. 

Let $M$ be a nonzero finitely generated graded  $R$-module and let
$$0 \rightarrow \bigoplus_{j\in\Z} R(-j)^{\beta_{p,j}(M)} \rightarrow \cdots \rightarrow \bigoplus_{j\in\Z}R(-j)^{\beta_{0,j}(M)}\rightarrow 0$$
be the minimal free resolution of $M$. The \emph{Castelnuovo–Mumford regularity} (or regularity for short) of $M$ is defined by
$$\reg(M) = \max\{j-i\mid \beta_{i,j}(M)\ne 0\},$$
and the {\it projective dimension} of $M$ is the length of this resolution
$$\pd(M) = p.$$

Let us denote by $d(M)$ the maximal degree of a minimal homogeneous generator of $M$. The definition of the regularity implies $$d(M) \leqslant \reg(M).$$

For any nonzero proper homogeneous ideal $I$ of $R$, by looking at the minimal free resolution, it is easy to see that $\reg(I) = \reg(R/I)+ 1$, so we shall work with $\reg(I)$ and $\reg(R/I)$ interchangeably.

The regularity of $M$ can also be computed via the local cohomology modules of $M$. For $i=0,\ldots, \dim(M)$, we define the $a_i$-invariant of $M$ as follows
$$a_i(M) = \max\{t \mid H_{\mi}^i(M)_t \ne 0\}$$
where $H_{\mi}^i(M)$ is the $i$-th local cohomology module of $M$ with the support $\mi = (x_1,\ldots,x_r)$ (with the convention $\max \emptyset = -\infty$). Then,
$$\reg(M) = \max\{a_i(M) + i\mid i = 0,\ldots, \dim(M)\},$$
and $$\pd(M) = r - \min\{i\mid H_{\mi}^i(M) \ne 0\}.$$

For example, since $\dim(R/\mi) = 0$ and $H_\mi^0(R/\mi) = R/\mi$, we have
$$\reg(\mi) = \reg(R/\mi)+1 = a_0(R/\mi) + 1 = \max\{i \mid \left(R/\mi\right)_i \ne 0\}+1 = 1.$$

\begin{rem} As usual we shall make the convention that $\reg(M) = -\infty$ if $M = 0$.\end{rem}

\subsection{Simplicial complexes and Stanley-Reisner ideals} A simplicial complex $\Delta$ over a finite set $V$ is a collection of subsets of $V$ such that if $F \in \Delta$ and $G\subseteq F$ then $G\in\Delta$. Elements of $\Delta$ are called faces. Maximal faces (with respect to inclusion) are called facets. For $F \in \Delta$, the dimension of $F$ is defined to be $\dim F = |F|-1$. The empty set, $\emptyset$, is the unique face of dimension $-1$, as long as $\Delta$ is not the void complex $\{\}$ consisting of no subsets of $V$. If every facet of $\Delta$ has the same cardinality, then $\Delta$ is called a {\it pure} complex. The dimension of $\Delta$ is $\dim \Delta = \max\{\dim F \mid F\in \Delta\}$. The link of $F$ inside $\Delta$ is its subcomplex:
$$\lk_{\Delta}(F) = \{H\in \Delta \mid H\cup F\in \Delta \ \text{ and } H\cap F=\emptyset\}.$$

Every element in a face of $\Delta$ is called a {\it vertex} of $\Delta$. Let us denote $V(\Delta)$ to be the set of vertices of $\Delta$. If there is a vertex, say $j$, such that $\{j\}\cup F \in \Delta$ for every $F\in \Delta$, then $\Delta$ is called a {\it cone} over $j$. It is well-known that if $\Delta$ is a cone, then it is an acyclic complex. A complex is called a {\it simplex} if it contains all subsets of its vertices, and thus a simplex is a cone over every its vertex.

For a subset $\tau = \{j_1,\ldots,j_i\}$ of $[r]$, denote $\x^\tau = x_{j_1} \cdots x_{j_i}$. Let $\Delta$ be a simplicial complex over the set $V=\{1,\ldots,r\}$. The Stanley-Reisner ideal of $\Delta$ is defined to be the squarefree monomial ideal
$$I_{\Delta} = (\x^\tau \mid \tau \subseteq [r] \text{ and } \tau \notin \Delta) \ \text{ in } R = K[x_1,\ldots,x_r]$$
and the {\it Stanley-Reisner} ring of $\Delta$ to be the quotient ring $k[\Delta] = R/I_{\Delta}$.  This provides a bridge between combinatorics and commutative algebra (see \cite{MS,ST}). 

Note that if $I$ is a square-free monomial ideal, then it is a Stanley-Reisner ideal of the simplicial complex $\Delta(I)= \{\tau \subseteq [r] \mid \x^\tau\not \in I\}$. When $I$ is a monomial ideal (maybe not square-free) we also use $\Delta(I)$ to denote the simplicial complex corresponding to the square-free monomial ideal $\sqrt{I}$. 

The regularity of a square-free monomial ideal can compute via the vanishing of reduced homology of simplicial complexes. From Hochster's formula on the Hilbert series of  the local cohomology module $H_\mi^i(I_\Delta)$ (see \cite[Corollary 13.16]{MS}), one has

\begin{lem}\label{Hochster-Reg} For a simplicial complex $\Delta$, we have
$$\reg(I_\Delta) = \max\{d\mid \h_{d-1}(\lk_{\Delta}(\sigma); K)\ne 0, \text{ for some } \sigma\in\Delta\}.$$
\end{lem}

The {\it Alexander dual} of $\Delta$, denoted by $\Delta^*$, is the simplicial complex over $V$ with faces
$$\Delta^* = \{V\setminus \tau\mid \ \tau\notin \Delta\}.$$
Notice that $(\Delta^*)^* = \Delta$. If $I = I_\Delta$ then we shall denote the Stanley-Reisner ideal of the Alexander dual $\Delta^*$ by $I^*$. It is a well-known result of Terai \cite{Terai} (or see \cite[Theorem 5.59]{MS})  that the regularity of a squarefree monomial ideal can be related to the projective dimension of its Alexander dual.
\begin{lem}\label{TeraiF} Let $I \subseteq R$ be a square-free monomial ideal. Then,
$$\reg(I) = \pd(R/I^*).$$
\end{lem}

Let $\F(\Delta)$ denote the set of all facets of $\Delta$.  We say that $\Delta$ is generated by $\F(\Delta)$ and write $\Delta = \left< \F(\Delta)\right>$. Note that $I_\Delta$ has the minimal primary decomposition (see \cite[Theorem 1.7]{MS}):
$$I_\Delta = \bigcap_{F\in \F(\Delta)} (x_i\mid i\notin F),$$
and therefore the $n$-th symbolic power of $I_\Delta$ is
$$I_\Delta^{(n)} = \bigcap_{F\in \F(\Delta)} (x_i\mid i\notin F)^n.$$

We next describe a formula to compute the local cohomology modules of monomial ideals.  Let $I$ be a non-zero monomial ideal. Since $R/I$ is an $\N^r$-graded algebra, $H^i_{\frak m}(R/I)$ is an $\mathbb Z^r$-graded module over $R/I$ for every $i$. For each degree $\alb=(\alpha_1,\ldots,\alpha_r)\in\Z^r$, in order to compute $\dim_K H_{\mi}^i(R/I)_{\alb}$ we use a formula given by Takayama \cite[Theorem $2.2$]{T} which is a generalization of Hochster's formula for the case $I$ is square-free \cite[Theorem 4.1]{ST}.

Set $G_{\alb}=\{i\mid \alpha_i<0\}$. For a subset $F\subseteq [r]$, we set $R_F=R[x_i^{-1}\mid i\in F\cup G_{\alb}]$. Define the simplicial complex $\Delta_{\alb}(I)$ by
\begin{equation}\label{degree-complex}
\Delta _{\alpha }(I) =\{F\subseteq [r]\setminus G_{\alb}\mid x^{\alpha }\notin IR_F\}.
\end{equation}

\begin{lem} \label{TA}\cite[Theorem 2.2]{T}  $\dim_K {H_{\frak m}^i(R/I)_{\alb}}=\dim_K \widetilde{H}_{i-\mid G_{\alpha }\mid-1 }(\Delta _{\alb}(I);K).$ 
\end{lem}

The following result of Minh and Trung is very useful to compute $\Delta_{\alb}(I_\Delta^{(n)})$, which allows us to investigate $\reg(I_\Delta^{(n)})$ by using the theory of convex polyhedra.

\begin{lem}\cite[Lemma 1.3]{MT} \label{MTr} Let $\Delta$ be a simplicial complex and $\alb \in\N^r$. Then,
$$\F(\Delta_{\alb}(I_\Delta^{(n)})) =\left\{F\in \F(\Delta) \mid \sum_{i\notin F} \alpha_i \leqslant n-1\right\}.$$ 
\end{lem}

This lemma can be generalized a little bit as follows.

\begin{lem}\cite[Lemma 1.3]{HT2} \label{HoaTrLem} Let $\Delta$ be a simplicial complex and $\alb \in\Z^r$. Then,
$$\F(\Delta_{\alb}(I_\Delta^{(n)})) =\left\{F\in \F(\lk_\Delta(G_{\alb})) \mid \sum_{i\notin F \cup G_{\alb}} \alpha_i \leqslant n-1\right\}.$$ 
\end{lem}

\subsection{Hypergraphs} \label{hyperSub} Let  $V$ be a finite set. A simple hypergraph $\H$ with vertex set $V$ consists of a set of subsets of $V$, called the edges of $\H$, with the property that no edge contains another. We use the symbols $V(\H)$ and $E(\H)$ to denote the vertex set and the edge set of $\H$, respectively. 

In this paper we assume that all hypergraphs are simple unless otherwise specified.

In the hypergraph $\H$, an edge is {\it trivial} if it contains only one element, a vertex is {\it isolated} if it is not appearing in any edge,  a vertex is a {\it neighbor} of another one if they are in some edge. 

A hypergraph $\H'$ is a {\it subhypergraph} of $\H$ if $V(\H') \subseteq V(\H)$ and $E(\H') \subseteq E(\H)$. For an edge $e$ of $\H$, we define $\H\setminus e$ to be the hypergraph obtained by deleting $e$ from the edge set of $\H$. For a subset $S\subseteq V(\H)$, we define $\H\setminus S$ to be the hypergraph obtained from $\H$ by deleting the vertices in $S$ and all edges containing any of those vertices.

A set $S\subseteq E(\H)$ is called an {\it edgewise dominant set} of $\H$ if every non-isolated vertex of $\H$ not contained in some edge of $S$ or contained in a trivial edge has a neighbor contained in some edge of $S$. Define,
$$\epsilon(\H) = \min\{|S| \mid S \text{ is edgewise dominant}\}.$$

For a hypergraph $\H$ with $V(\H)\subseteq [r]$, we associate to the hypergraph $\H$ a square-free monomial ideal
$$I(\H) = (\x^e\mid e\in E(\H)) \subseteq R,$$
which is called the {\it edge ideal} of $\H$.

Notice that if $I$ is a square-free monomial ideal, then $I$ is an edge ideal of a hypergraph with the edge set uniquely determined  by the generators of $I$.

Let $\H^*$ be the simple hypergraph corresponding to the Alexander duality $I(\H)^*$ of $I(\H)$. We will determine the edge set of $\H^*$, it turns out that $E(\H^*)$ is the set of all minimal vertex covers of $\H$. A {\it vertex cover} in a hypergraph is a set of vertices, such that every edge of the hypergraph contains at least one vertex of that set. It is an extension of the notion of vertex cover in a graph.  A vertex cover $S$ is called minimal if no proper subset of  $S$ is a vertex cover. From the minimal primary decomposition (see \cite[Definition 1.35 and Proposition 1.37]{MS}):
$$I(\H^*) = \bigcap_{e\in E(\H)} (x_i\mid i\in e),$$
it follows that $E(\H^*)$ is just the set of minimal vertex covers of $\H$. Thus,
$$I(\H^*) = (\x^\tau\mid \tau \text{ is a minimal vertex cover of } \H).$$


In the sequel, we need the following result of Dao and Schweig \cite[Theorem 3.2]{DS}.

\begin{lem}\label{DS} Let $\H$ be a hypergraph. Then, $\pd(R/I(\H)) \leqslant |V(\H)| - \epsilon(\H)$.
\end{lem}

\subsection{Matchings in a graph} Let $G$ be a graph. A {\it matching} in $G$ is a subgraph consisting of pairwise disjoint edges. If this subgraph is an induced subgraph, then the matching is called an {\it induced matching}. A matching of $G$ is maximal if it is maximal with respect to inclusion. The {\it matching number} of G, denoted by $\match(G)$, is the maximum size of a matching in $G$; and the {\it induced matching number} of $G$, denoted by $\nu(G)$, is the maximum size of an induced matching in $G$.

An {\it independent set} in $G$ is a set of vertices no two of which are adjacent to each other. An independent set in $G$ is maximal (with respect to set inclusion) if the set cannot be extended to a larger independent set. Let $\Delta(G)$ denote the set of all independent sets of $G$. Then, $\Delta(G)$ is a simplicial complex, called the {\it independence complex} of $G$. It is well-known that $I(G) = I_{\Delta(G)}$.

According to Constantinescu and Varbaro \cite{CV}, we say that a matching $M=\{\{u_i,v_i\} \mid i=1,\ldots,s\}$ is an {\it ordered matching} if:
\begin{enumerate}
\item $\{u_1,\ldots,u_s\}\in \Delta(G)$,
\item $\{u_i, v_j\} \in E(G)$ implies $i \leqslant j$.
\end{enumerate}
The {\it ordered matching number} of $G$, denoted by $\ordmatch(G)$ is the maximum size of an ordered matching in $G$.

\medskip

The following result gives a lower bound for $\reg(I(G)^{(n)})$ in terms of the induced matching number $\nu(G)$ 

\begin{lem}\label{ind-reg}\cite[Theorem 4.6]{GHOS} Let $G$ be a graph. Then,
$$\reg(I(G)^{(n)}) \geqslant 2n + \nu(G)-1, \text{ for all } n\geqslant 1.$$
\end{lem}

\subsection{Convex polyhedra} The theory of convex polyhedra plays a key role in our study.

For a vector $\alb = (\alpha_1,\ldots,\alpha_r)\in\R^r$, we set $|\alb| : = \alpha_1+\cdots+\alpha_r$ and for a nonempty bounded closed subset $S$ of $\R^r$ we set 
$$\delta(S) : = \max\{|\alb| \mid \alb \in S\}.$$

Let $\Delta$ be a simplicial conplex over $[r]$. In general, $\reg(I_\Delta^{(n)})$ is not a linear function in $n$ for $n\gg 0$ (see e.g. \cite[Theorem 5.15]{DHHT}), but a quasi-linear function as in the following result.

\begin{lem}\label{HT-bound} \cite[Theorem 4.9]{HT1} There exist positive integers $N, n_0$ and rational numbers $a,b_0,\ldots,b_{N-1} < \dim(R/I_\Delta)+1$ such that 
$$\reg(I_\Delta^{(n)})=an+b_k, \text{ for all } n \geqslant n_0 \text{ and } n\equiv k \mod N,  \text{ where } 0\leqslant k \leqslant N-1.$$
Moreover, $\reg(I_\Delta^{(n)}) < an + \dim(R/I_\Delta)+1$ for all $n\geqslant 1$.
\end{lem}

By virtue of this result, we define
$$\delta(I_\Delta) = a = \lim\limits_{n\to\infty} \frac{\reg(I_\Delta^{(n)})}{n}.$$

In order to compute this invariant we can use the geometric interpretation of it by means of  symbolic polyhedra defined in \cite{CEHH, DHHT}. Let $\smp(I_{\Delta})$ be the convex polyhedron in $\R^r$ defined by the following system of linear inequalities:
\begin{equation}\label{EQ1}
\begin{cases}
\sum\limits_{i\notin F} x_i \geqslant 1 & \text{for}~ F\in\mathcal F(\Delta),\\
x_1\geqslant 0,\ldots,x_r\geqslant 0,
\end{cases}
\end{equation}
which is called the {\it symbolic polyhedron} of $I_\Delta$. Then, $\smp(I_{\Delta})$ is a convex polyhedron in $\R^r$. By \cite[Theorem 3.6]{DHHT} we have 
\begin{equation}\label{lead-sym-ideal}
\delta(I_\Delta) = \max\{|\v| \mid \v \text{ is a vertex of } \smp(I_{\Delta})\}.
\end{equation}

Now assume that  
$$H_\mi^i(I_\Delta^{(n)})_{\alb} \ne 0$$
for some $0\leqslant i \leqslant \dim(R/I_{\Delta})$ and $\alb = (\alpha_1,\ldots,\alpha_r) \in \N^r$.

By Lemma \ref{TA} we have
\begin{equation}\label{SH1}
\dim_K \h_{i-1}(\Delta_{\alb}(I_\Delta^{(n)});K) =\dim_K H_{\mi}^i(R/I_\Delta^{(n)})_{\alb} \ne 0.
\end{equation}
In particular, $\Delta_{\alb}(I_\Delta^{(n)})$ is not acyclic.

Suppose that $\F(\Delta) = \{F_1,\ldots,F_t\}$ for $t\geqslant 1$. By Lemma \ref{MTr} we may assume that
$$\F(\Delta_{\alb}(I_\Delta^{(n)})) = \{F_1,\ldots,F_s\}, \text{ where }1\leqslant s \leqslant t.$$

For each integer $m\geqslant 1$, let $\P_m$ be the convex polyhedron of $\R^r$ defined by:
\begin{equation}\label{EQ-basics}
\begin{cases}
\sum\limits_{i\notin F_j} x_i  \leqslant m-1 & \text{ for } j = 1,\ldots,s,\\
\sum\limits_{i\notin F_j} x_i  \geqslant  m & \text{ for } j = s+1,\ldots,t,\\
x_1\geqslant 0,\ldots,x_r\geqslant 0.
\end{cases}
\end{equation}
Then, $\alb \in \mathcal P_n$. Moreover, by Lemma $\ref{MTr}$ one has
\begin{equation} \label{SH2}
\Delta_{\btb}(I_\Delta^{(m)}) = \left<F_1,\ldots,F_s\right> = \Delta_{\alb}(I_\Delta^{(n)}) \ \text{ whenever } \btb\in \P_m \cap \N^r.
\end{equation}
Note also that for such a vector $\btb$, by Formula $(\ref{SH2})$ we have
$$\dim_K \h_{i-1}(\Delta_{\btb}(I_\Delta^{(m)});K) = \dim_K \h_{i-1}(\Delta_{\alb}(I_\Delta^{(n)});K) \ne 0.$$
Together with Lemma \ref{TA}, this fact yields
\begin{equation} \label{SH3}
H_\mi^i(R/ I_\Delta^{(m)})_{\btb} \ne 0.
\end{equation}

In order to investigate the convex polyhedron $\P_m$ we also consider the convex polyhedron $\C_m$ in $\R^r$ defined by:
\begin{equation}\label{EQ-polytope}
\begin{cases}
\sum\limits_{i\notin F_j} x_i  \leqslant m & \text{ for } j = 1,\ldots,s,\\
\sum\limits_{i\notin F_j} x_i  \geqslant  m & \text{ for } j = s+1,\ldots,t,\\
x_1\geqslant 0,\ldots,x_r\geqslant 0.
\end{cases}
\end{equation}
Note that $\C_m = m\C_1$ for all $m\geqslant 1$,  where $m\C_1= \{m\y \mid \y \in \C_1\}$. 

By the same way as in the proof of \cite[Lemma 2.1]{HT0} we obtain the following lemma.

\begin{lem}\label{polytope} $\C_1$ is a polytope with $\dim \C_1 = r$.
\end{lem}

The next lemma gives an upper bound for $\delta(\C_1)$.

\begin{lem}\label{bound-lead} $\delta(C_1) \leqslant \delta(I_\Delta)$.
\end{lem}
\begin{proof} Since $\C_1$ is a polytope with $\dim \C_1 = r$ by Lemma \ref{polytope}, $\delta(\C_1) = |\gmb|$ for some vertex $\gmb$ of $\C_1$. By \cite[Formula $(23)$ in Page 104]{S} we imply that $\gmb$ is the unique solution of a system of linear equations of the form
\begin{equation}\label{EQ-gamma}
\begin{cases}
\sum\limits_{i\notin F_j} x_i  = 1 & \text{ for } j \in S_1,\\
x_j=0  & \text{ for } j \in S_2,
\end{cases}
\end{equation}
where $S_1 \subseteq [t]$ and $S_2\subseteq [r]$ such that $|S_1|+|S_2| = r$. By using Cramer's rule to get $\gmb$, we conclude that $\gmb$ is a rational vector. In particular, there is a positive integer, say $p$, such that $p\gmb \in \N^r$. Note that $\C_p = p \C_1$, so $p\gmb \in \C_p \cap \N^r$.

For every $j\geqslant 1$, let $\y = jp\gmb + \alb$. Then, $\y \in \N^r$ and $|\y|= \delta(\C_1)jp + |\alb|$. On the other hand, by using the fact that $jp\gmb \in \C_{jp}$, we can check that
$$
\begin{cases}
\sum\limits_{i\notin F_j} y_i  \leqslant jp+n-1 & \text{ for } j = 1,\ldots,s,\\
\sum\limits_{i\notin F_j} y_i  \geqslant  jp+n & \text{ for } j = s+1,\ldots,t,\\
\end{cases}
$$
and so $\y \in \P_{jp+n} \cap \N^r$. 

Together with Equation $(\ref{SH3})$, we deduce that  $H_\mi^i(R/I_\Delta^{(jp+n)})_{\y} \ne 0$, and therefore
$$\reg(R/I_\Delta^{(jp+n)}) \geqslant |\y| + i =  \delta(\C_1)jp + |\alb| + i.$$
Combining with Lemma \ref{HT-bound}, this inequality yields
$$\delta(\C_1)jp + |\alb| + i <  \delta(I_\Delta)(jp+n) + \dim(R/I_\Delta).$$
Since this inequality valid for any positive integer $j$, it forces $\delta(\C_1) \leqslant \delta(I_\Delta)$.
\end{proof}

\section{Regularity of symbolic powers of  ideals}

In this section we will prove the upper bound for $\reg(I_\Delta^{(n)})$.  Firts we start with the following fact.

\begin{lem} \label{complex-one} Let $\sigma \subseteq [r]$ with $\sigma\ne [r]$, $S = K[x_i\mid i\notin \sigma]$ and $J = IR_\sigma \cap S$. Then,
$$\reg(J^{(n)}) \leqslant \reg (I^{(n)}) \ \text{ for all } n\geqslant 1.$$ 
In particular, $\delta(J) \leqslant \delta(I)$.
\end{lem}
\begin{proof} We may assume that $S = K[x_1,\ldots, x_s]$ for some $1\leqslant s \leqslant r$. Let $i$ be an index and $\alb$ a vector in $\Z^s$ such that
$$H_{\ni}^i(S/J^{(n)}) _{\alb}\ne 0 \text{ and }  \reg(S/J^{(n)}) = |\alb| + i,$$ 
where $\ni = (x_1,\ldots,x_s)$ is the homogeneous maximal ideal of $S$. 

Let $\btb = (\alpha_1,\ldots,\alpha_s, -1,\ldots,-1)\in \Z^r$ so that $G_{\btb} = G_{\alb} \cup\{s+1,\ldots,r\}$. By Formula $(\ref{degree-complex})$ we deduce that
\begin{equation}\label{EQLM1}
\Delta_{\alb}(J^{(n)}) = \Delta_{\btb}(I^{(n)}).
\end{equation}

By Lemma \ref{TA}, $$\dim_K H_{\ni}^i(S/J^{(n)}) _{\alb} = \dim_K \h_{i-|G_{\alb}|-1}(\Delta_{\alb}(J^{(n)});K),$$
and thus $\h_{i-|G_{\alb}|-1}(\Delta_{\alb}(J^{(n)});K) \ne 0$. Together with Equation (\ref{EQLM1}), it yields
$$\h_{i-|G_{\alb}|-1}(\Delta_{\btb}(I^{(n)});K) \ne 0.$$

By Lemma \ref{TA} again, it gives $H_{\mi}^{i+(r-s)}(R/I^{(n)}) _{\btb}\ne 0$ since $|G_{\btb}| = |G_{\alb}| + (r-s)$. Therefore,
$$\reg(R/I^{(n)}) \geqslant |\btb| + i+(r-s) = |\alb| + i = \reg(S/J^{(n)}),$$
it follows that $\reg(J^{(n)}) \leqslant \reg(I^{(n)})$.

Finally, together this inequality with Lemma \ref{HT-bound} we have
$$\delta(J) = \lim\limits_{n\to\infty} \frac{\reg(J^{(n)})}{n} \leqslant \lim\limits_{n\to\infty} \frac{\reg(I^{(n)})}{n} = \delta(I),$$
and the lemma follows.
\end{proof}

\begin{thm}\label{a-Invariant} Let $I$ be a square-free monomial ideal. Then, for all $i \geqslant 0$ we have
$$a_i(R/I^{(n)}) \leqslant \delta(I)(n -1).$$
\end{thm}
\begin{proof} If $n = 1$, the theorem follows from Hochster's formula on the Hilbert series of  the local cohomology module $H_\mi^i(R/I_\Delta)$  (see \cite[Theorem 4.1]{ST}).

We may assume that $n\geqslant 2$. If $a_i(R/I^{(n)}) = -\infty$, the theorem is obvious, so that we also assume that $a_i(R/I^{(n)}) \ne -\infty$.

Suppose  $\alb \in \Z^r$ such that
$$H_{\mi}^i(R/I^{(n)}) _{\alb}\ne 0 \text{ and }  a_i(R/I^{(n)}) = |\alb|.$$ 

By Lemma $\ref{TA}$ we have
\begin{equation}\label{N1}
\dim_K \h_{i-|G_{\alb}|-1}(\Delta_{\alb}(I^{(n)});K) =\dim_K H_{\mi}^i(R/I^{(n)})_{\alb} \ne 0.
\end{equation}
In particular, $\Delta_{\alb}(I^{(n)})$ is not acyclic.

If $G_{\alb} =[r]$, then $a_i(R/I^{(n)}) = |\alb|\leqslant 0$, and hence the theorem holds in this case.

We therefore assume that $G_{\alb} =\{m+1,\ldots,r\}$ for $1\leqslant m\leqslant  r$. Let $S = K[x_1,\ldots,x_m]$ and  $J = IR_{G_{\alb}} \cap S$.

Let $\alb' = (\alpha_1,\ldots,\alpha_m) \in \N^m$. By using Formula $(\ref{degree-complex})$, we have
\begin{equation}\label{POWER}
\Delta_{\alb'}(J^{(n)}) = \Delta_{\alb}(I^{(n)}).
\end{equation}
Together with $(\ref{N1})$, it gives $\h_{i-|G_{\alb}|-1}(\Delta_{\alb'}(J^{(n)});K)\ne 0$. By Lemma $\ref{TA}$ we get
$$H_{\ni}^{i-|G_{\alb}|}(S/J^{(n)})_{\alb'}\ne 0,$$
where $\ni = (x_1,\ldots,x_m)$ is the homogeneous maximal ideal of $S$. 

Let $\Delta$ be the simplicial complex over $[m]$ corresponding to the square-free monomial ideal $J$. Assume that $\F(\Delta) = \{F_1,\ldots,F_t\}$.

By Lemma \ref{MTr} we may assume that $\F(\Delta_{\alb'}(J^{(n)})) = \{F_1,\ldots,F_s\}$ for $1\leqslant s \leqslant t$. Let $$\btb =(\beta_1,\ldots,\beta_m)= \dfrac{1}{n-1}\alb' \in \R^m.$$

By Lemma \ref{MTr} again, we deduce that
$$
\begin{cases}
\sum\limits_{i\not \in F_j} \beta_i =\dfrac{1}{n-1} \sum\limits_{i\not \in F_j} \alpha_i  \leqslant 1 & \text{ for } j = 1,\ldots,s,\\
\sum\limits_{i\not \in F_j} \beta_i =\dfrac{1}{n-1} \sum\limits_{i\not \in F_j} \alpha_i  \geqslant \dfrac{n}{n-1} > 1 & \text{ for } j = s+1,\ldots,t.\\
\end{cases}
$$
It follows that $\btb\in C_1$, where $C_1$ is a polyhedron in $\R^m$ defined by
$$
\begin{cases}
\sum\limits_{i\not \in F_j} x_i  \leqslant 1 & \text{ for } j = 1,\ldots,s,\\
\sum\limits_{i\not \in F_j} x_i  \geqslant  1 & \text{ for } j = s+1,\ldots,t,\\
x_1\geqslant 0,\ldots,x_m\geqslant 0.
\end{cases}
$$
By Lemma \ref{polytope}, $C_1$ is a polytope in $\R^m$.

Hence $|\btb| \leqslant \delta(C_1)$, and hence $|\alb'|  = (n-1)|\btb| \leqslant \delta(C_1)(n-1)$. Observe that $\alpha_j < 0$ for all $j \in G_{\alb} = \{m+1, \ldots,r\}$, so
\begin{equation}\label{MF1}
a_i(R/I^{(n)}) = |\alb| = |\alb'| + (\alpha_{m+1}+\cdots+\alpha_r) \leqslant |\alb'| \leqslant \delta(C_1)(n-1). 
\end{equation}

On the other hand, by Lemmas \ref{bound-lead} and \ref{complex-one} we deduce that 
$$\delta(C_1)\leqslant \delta(J) \leqslant \delta(I).$$

Together with Formula $(\ref{MF1})$, it yields $a_i(R/I^{(n)}) \leqslant \delta(I)(n-1)$, and the proof of the theorem is complete.
\end{proof}

We are now in position to prove the main result of the paper.

\begin{thm}\label{main-thm} Let $\Delta$ be a simplicial complex. Then,
$$\reg(I_\Delta^{(n)}) \leqslant \delta(I_\Delta)(n -1)+ b, \ \text{ for all } n\geqslant 1,$$
where 
$b = \max\{\reg(I_\Gamma) \mid  \Gamma \text{ is a subcomplex of } \Delta \text{ with }  \F(\Gamma) \subseteq \F(\Delta)\}$.
\end{thm}
\begin{proof} For simplicity, we put $I = I_{\Delta}$. Let $i \in\{0, \ldots, \dim(R/I)\}$ and $\alb\in\Z^r$ such that
$$H_\mi^i(R/I^{(n)})_{\alb} \ne 0, \text{ and } \reg(R/I^{(n)}) = a_i(R/I^{(n)})+i = |\alb| + i.$$

By Lemma \ref{TA}, we have
\begin{equation}\label{N10}
\dim_K \h_{i-|G_{\alb}|-1}(\Delta_{\alb}(I^{(n)});K) =\dim_K H_{\mi}^i(R/I^{(n)})_{\alb} \ne 0.
\end{equation}
In particular, $\Delta_{\alb}(I^{(n)})$ is not acyclic.

If $G_{\alb} =[r]$, then $\Delta_{\alb}(I^{(n)})$ is either $\{\emptyset\}$ or a void complex. Because it is not acyclic, $\Delta_{\alb}(I^{(n)}) =\{\emptyset\}$. By Formula $(\ref{N10})$ we deduce that $i = |G_{\alb}| = r$, and hence $\dim R/I = r$. It means that $I = 0$,  so $I^{(n)} = 0$ as well. Therefore,  $\reg(I^{(n)}) = -\infty$, and the theorem holds in this case.

We may assume that $G_{\alb} = \{m+1,\ldots,r\}$ for some $1\leqslant m\leqslant r$. Let $S = K[x_1,\ldots,x_m]$ and  $J = IR_{G_{\alb}} \cap S$.

 Let $\alb' = (\alpha_1,\ldots,\alpha_m) \in \N^m$. By using Formula $(\ref{degree-complex})$, we have
\begin{equation}\label{POWER02}
\Delta_{\alb'}(J^{(n)}) = \Delta_{\alb}(I^{(n)}).
\end{equation}
Together with $(\ref{N10})$, it gives $\h_{i-|G_{\alb}|-1}(\Delta_{\alb'}(J^{(n)});K)\ne 0$. By Lemma $\ref{TA}$ we get
$$H_{\ni}^{i-|G_{\alb}|}(S/J^{(n)})_{\alb'}\ne 0,$$
where $\ni = (x_1,\ldots,x_m)$ is the homogeneous maximal ideal of $S$. In particular,
$$|\alb'| \leqslant a_{i -|G_{\alb}|}(S/J^{(n)}).$$
Together with Lemma \ref{complex-one} and  Theorem \ref{a-Invariant}, it yields
$$|\alb'|  \leqslant \delta(J)(n-1) \leqslant \delta(I)(n-1).$$
Therefore,
$$\reg(I^{(n)}) = |\alb| + i =|\alb'| +\sum_{j=m+1}^r \alpha_j +i \leqslant |\alb'| + i - |G_{\alb}|\leqslant \delta(I)(n-1) + i - |G_{\alb}|.$$

It remains to prove that $i - |G_{\alb}| \leqslant b$. By  Lemma \ref{HoaTrLem},  we have
$$\Delta_{\alb'}(J^{(n)}) =\Delta_{\alb}(I^{(n)})= \left\{F\in\F(\lk_{\Delta}(G_{\alb}))\mid \sum_{j\notin F\cup G_{\alb}}\alb_j \leqslant n-1  \right\}.$$
It follows that there is a simplicial complex $\Gamma$ with $\F(\Gamma) \subseteq \F(\Delta)$ such that $$\Delta_{\alb'}(J^{(n)}) = \lk_{\Gamma}(G_{\alb}).$$

Since $\h_{i-|G_{\alb}|-1}(\lk_{\Gamma}(G_{\alb});K) \ne 0$, by Lemma \ref{Hochster-Reg} we have $i-|G_{\alb}| \leqslant \reg(I_{\Gamma}) \leqslant b$, and then proof of the theorem is complete.
\end{proof}

As a direct consequence of Theorem \ref{main-thm}, we have a simple bound. Namely,

\begin{cor}\label{main-cor} Let $I$ be a square-free monomial ideal. Then, 
$$\reg(I^{(n)}) \leqslant \delta(I)(n -1)+ \dim(R/I)+1, \ \text{ for all } n\geqslant 1.$$
\end{cor}
\begin{proof} Let $\Delta$ be the simplicial complex corresponding to the square-free ideal $I$. For every subcomplex $\Gamma$ of $\Delta$ we have $\dim\Gamma \leqslant \dim\Delta$. It follows from Lemma \ref{Hochster-Reg} that
$$\reg(I_{\Gamma}) \leqslant \dim (R/I_{\Gamma}) + 1 \leqslant \dim(R/I_\Delta)+1.$$
Therefore, $b = \max\{\reg(I_\Gamma) \mid \F(\Gamma) \subseteq \F(\Delta)\} \leqslant \dim(R/I_\Delta)+1$. Now the corollary follows from Theorem \ref{main-thm}.
\end{proof}

\medskip

We next reformulate the theorem \ref{main-thm} for a square-free monomial ideal arising from a hypergraph.

\begin{thm}\label{hyper-state} Let $\H$ be a hypergraph. Then, for all $n\geqslant 1$, we have
$$\reg(I(\H)^{(n)}) \leqslant \delta(I(\H))(n-1) + b,$$
where $b = \max\{\pd(R/I(\H'))\mid \H' \text{ is a subhypergraph of } \H^* \text{ with } E(\H') \subseteq E(\H^*)\}$.
\end{thm}
\begin{proof} Let $\Delta$ be the corresponding simplicial complex of the square-free monomial ideal $I(\H)$. Assume that $\F(\Delta) = \left\{F_1,\ldots,F_p\right\}$. Since
$$I(\H) = \bigcap_{j=1}^p  (x_i \mid i \notin F_j),$$
so that $E(\H^*) = \{C_1,\ldots, C_p\}$, where $C_j = [r] \setminus F_j$ for all $j=1,\ldots,p$.

Let $\Gamma$ be a subcomplex of $\Delta$ with $\F(\Gamma) \subseteq \F(\Delta)$. We may assume that $\F(\Gamma) = \{F_1,\ldots,F_k\}$ for $1\leqslant k \leqslant p$. Then, we have
$I_\Gamma^* = I(\H')$ where $\H'$ is the subhypergraph of $\H^*$ with $E(\H') = \{C_1,\ldots,C_k\}$.

By Lemma \ref{TeraiF} we have $\reg(I_\Gamma) = \pd(R/I_\Gamma^*) = \pd(R/I(\H'))$, and therefore the theorem follows from Theorem \ref{main-thm}.
\end{proof}

The next theorem is the second main result of the paper. It bounds  the regularity of symbolic powers of a square-free monomial ideal via the combinatorial properties of the associated hypergraph.

\begin{thm}\label{main-hyper} Let $\H$ be a simple hypergraph. Then,
$$\reg(I(\H)^{(n)}) \leqslant \delta(I(\H))(n-1) + |V(\H)| - \epsilon(\H^*), \text{ for all } n\geqslant 1.$$
\end{thm}
\begin{proof} By Theorem \ref{hyper-state}, it suffices to show that
$$\pd(R/I(\G)) \leqslant |V(\H)| - \epsilon(\H^*)$$
for every hypergraph $\G$ with $E(\G) \subseteq E(\H^*)$. By Lemma \ref{DS}, it suffices to prove that
$$|V(\G)| -\epsilon(\G) \leqslant |V(\H^*)|-\epsilon(\H^*).$$

In order to prove this inequality, without loss of generality we may assume that $\H^*$ has no both trivial edges and isolated vertices. 

Let $S$ be an edgewise-dominant set of $\G$ such that $|S| = \epsilon(\G)$. For each vertex $v\in V(\H^*)\setminus V(\G)$, we take an edge of $\H^*$ containing $v$, and denote this edge by $F(v)$. Then,
$$S' = S \cup \{F(v)\mid v\in V(\H^*) \setminus V(\G)\}$$
is an edgewise-dominant set of $\H^*$. It follows that 
$$\epsilon(\H^*) \leqslant |S'| \leqslant  |S| + |V(\H^*) \setminus V(\G)| = |S| + |V(\H^*)| - |V(\G)|,$$
and therefore $|V(\G)|-\epsilon(\G) \leqslant|V(\H^*)|-\epsilon(\H^*)$, as required.
\end{proof}

The following example shows that the bound in Theorem \ref{main-thm} is sharp at every $n$ for the class of matroid complexes. Recall that a simplicial complex $\Delta$ is called a {\it matroid complex} if for every subset $\sigma$ of $V(\Delta)$, the simplicial complex $\Delta[\sigma]$ is pure (see e.g. \cite[Chapter 3]{ST}). Here, $\Delta[\sigma]$ is the restriction of $\Delta$ to $\sigma$ and defined by $\Delta[\sigma] = \{\tau \mid \tau\in \Delta \text{ and } \tau \subseteq \sigma\}$.

\begin{exm}\label{exam-matroid} Let $\Delta$ be a matroid complex that is not a cone. Then,
$$\reg(I_\Delta^{(n)}) = \delta(I_\Delta)(n -1)+ b, \ \text{ for all } n\geqslant 1,$$
where 
$b = \max\{\reg(I_\Gamma) \mid  \Gamma \text{ is a subcomplex of } \Delta \text{ with }  \F(\Gamma) \subseteq \F(\Delta)\}$.
\end{exm}
\begin{proof} Let $I = I_\Delta$ and $s = \dim(R/I_\Delta)$. By \cite[Theorem 4.5]{MT2}, for all $n\geqslant 1$ we have:
$$\reg(I^{(n)}) = d(I)(n-1) +s+1.$$
It implies that
$$\lim\limits_{n\to\infty} \frac{\reg(I^{(n)})}{n}= d(I),$$
so $\delta(I) = d(I)$. It remains to show that $b = s+1$.

Together the fact $\delta(I) = d(I)$ with Theorem $\ref{main-thm}$, we get $s+1 \leqslant b$. On the other hand, by the same argument as in the proof of Corollary \ref{main-cor}, we obtain $b\leqslant s + 1$. Hence, $b = s + 1$, as required.
\end{proof}

We conclude this section with a remark on lower bounds.

\begin{rem} Let $I$ be a square-free monomial ideal. By \cite[Lemma 4.2(ii)]{DHHT} we deduce that $d(I) n \leqslant d(I^{(n)})$, and therefore
$$\reg (I^{(n)}) \geqslant d(I)n, \ \text{ for all } n\geqslant 1.$$
In general, $d(I) < \delta(I)$ (see e.g. \cite[Lemma 5.14]{DHHT}), so that the bound is not optimal. 

On the other hand, by Lemma \ref{HT-bound}, there is a number $b$ such that
$$\reg(I^{(n)}) \geqslant \delta(I)n + b, \text{ for all } n\geqslant 1.$$
The natural question is to find a good bound for $b$.
\end{rem}

\section{Applications}

In this section we will apply Theorem \ref{main-thm} to the regularity of symbolic powers of the edge ideal of a graph.  We start with a result which allows us to bound the number $b$ in Theorem \ref{main-thm} by choosing a suitable numerical function, it is of independent interest.

\begin{thm}\label{thm-delta} Let $\Delta$ be a simplicial complex over $[r]$ and let $$\simp(\Delta)=\{ \lk_{\Delta}(\sigma) \mid \sigma \in \Delta\}.$$
 Assume that $f\colon \simp(\Delta) \to \N$ is a function which satisfies the following properties:
\begin{enumerate}
\item If $\Lambda\in\simp(\Delta)$ is a simplex, then $f(\Lambda) = 0$.
\item For every $\Lambda\in \simp(\Delta)$ and every $v\in V(\Lambda)$ such that $\Lambda$ is not a cone over $v$, $f(\lk_\Lambda(v)) +1 \leqslant f(\Lambda)$.
\end{enumerate}
Then, for every subcomplex $\Gamma$ of $\Delta$ with $\F(\Gamma) \subseteq \F(\Delta)$ we have $\reg(I_\Gamma) \leqslant f(\Delta)+1$.
\end{thm}

\begin{proof} For a subset $S$ of $[r]$ we set $\mathfrak p_S = (x_i \mid i\in S) \subseteq R$. In order to facilitate an induction argument on the number of vertices of $\Delta$ we prove the following assertion:
\begin{equation}\label{main-graph}
\reg(\mathfrak{p}_S + I_\Gamma) \leqslant f(\Delta)+1, \ \text{ for every } S\subseteq [r],
\end{equation}
where all simplicial complexes is considered over $[r]$.

Indeed, if $|V(\Delta)| \leqslant 1$, then $\Delta$ is a simplex. In this case, the assertion is obvious.

Assume that $|V(\Delta)|\geqslant 2$. If $\Delta$ is a simplex, the assertion holds, so we assume that $\Delta$ is not a simplex. We now prove by backward induction  on $|S|$. If $|S| = r$, then
$$\mathfrak{p}_S + I_\Gamma = (x_1,\ldots,x_r).$$
In this case $\reg(\mathfrak{p}_S + I_\Gamma) = 1$, and so the assertion holds.

Assume that $|S| < r$. If $\mathfrak{p}_S + I_\Gamma$ is a prime, i.e. it is generated by variables, then $\reg(\mathfrak{p}_S + I_\Gamma) = 1$, and then the assertion holds.

Assume that  $\mathfrak{p}_S + I_\Gamma$ is not a prime. Then, there is  a variable, say $x_v$ with $v\in [r]$, such that $x_v$ appears in some monomial generator of $\mathfrak{p}_S + I_\Gamma$ of order at least $2$ and $v\notin S$. Note that if $u$ is not a vertex of $\Gamma$ then $x_u$ is a monomial generator of $I_\Gamma$, and if $\Gamma$ is a cone over some vertex $w$ then $x_w$ does not appear in any monomial generator of $I_\Gamma$. It implies that $v$ is a vertex of $\Gamma$ and $\Gamma$ is not a cone over $v$. In particular, $\Delta$ is not a cone over $v$.

Since
$$(\mathfrak{p}_S + I_\Gamma) + (x_v) = \mathfrak{p}_{S\cup \{v\}} + I_\Gamma, \text{ and } (\mathfrak{p}_S + I_\Gamma): (x_v) = \mathfrak{p}_S + I_{\Gamma'},$$
where $\Gamma'$ is a subcomplex of $\Gamma$ with $\F(\Gamma') = \{F \in  \F(\Gamma) \mid v \in F\}$, by \cite[Lemma 2.10]{DHS} we have
\begin{equation}\label{B001}
\reg(\mathfrak{p}_S + I_\Gamma) \leqslant \max\{\reg(\mathfrak{p}_{S\cup\{v\}} + I_\Gamma), \reg(\mathfrak{p}_S +  I_{\Gamma'})+1\}.
\end{equation}

By the backward induction hypothesis, we have  
\begin{equation}\label{B002}
\reg(\mathfrak{p}_{S\cup\{v\}} + I_\Gamma)\leqslant f(\Delta)+1.
\end{equation}

We now claim that
\begin{equation}\label{B003}
\reg(\mathfrak{p}_S +  I_{\Gamma'})\leqslant f(\Delta).
\end{equation}
Indeed, if $\mathfrak{p}_S +  I_{\Gamma'}$ is prime, then $\reg(\mathfrak{p}_S +  I_{\Gamma'})=1$. As $\Delta$ is not a cone over $v$, by the definition of $f$ we have
$f(\Delta) \geqslant f(\lk_\Delta(v))+1\geqslant 1$, and the claim holds in this case. 

Assume that $\mathfrak{p}_S +  I_{\Gamma'}$ is not a prime. Observe that $$I_{\Gamma''} = (x_v) + I_{\Gamma'},$$
where $\Gamma'' = \lk_{\Gamma'}(v)$ and this simplicial complex is considered over $[r]$. Since variable $x_v$ does not appear in any generator of $I_{\Gamma'}$, hence $\reg(I_{\Gamma''}) = \reg(I_{\Gamma'})$. 

On the other hand, by the induction hypothesis, we have 
$$\reg(I_{\Gamma''}) =\reg(\lk_{\Gamma'}(v)) \leqslant f(\lk_{\Delta}(v))+1.$$ 
It follows that
$$\reg(\mathfrak{p}_S + I_{\Gamma'}) \leqslant \reg(I_{\Gamma'})  =\reg(I_{\Gamma''}) \leqslant f(\lk_{\Delta}(v))+1.$$
Together with the inequality $f(\lk_{\Delta}(v))+1 \leqslant f(\Delta)$, it yields $\reg(\mathfrak{p}_S + I_{\Gamma'}) \leqslant f(\Delta)$, as claimed.

By combining three Inequalities (\ref{B001})-(\ref{B003}), we obtain $\reg(\mathfrak{p}_S + I_\Gamma) \leqslant f(\Delta)+1$, and so the inequality (\ref{main-graph}) is proved. The lemma now follows from the assertion by taking $S = \emptyset$, and the proof is complete.
\end{proof}

We now reformulate the theorem \ref{thm-delta} for graphs. A graph $G$ is called {\it trivial} if it has no edges. For a subset $S$ of $V(G)$, the {\it closed neighborhood} of the set $S$ in $G$ is the set $N_G[S] = S \cup \{v\in V(G) \mid v \text{ is a neighbor of some vertex in } S\}$. For a vertex $v$ of $G$,  we write $N_G[v]$ stands for $N_G[\{v\}]$. Recall that $\Delta(G)$ is the set of independent sets of $G$, which is a simplicial complex and $I(G) = I_{\Delta(G)}$.

\begin{cor}\label{lem-f} Let $G$ be a graph and let $\mathcal I_G=\{G\setminus N_G[S] \mid S\in\Delta(G)\}$. Assume that $f\colon \mathcal I_G \to \N$ is a function which satisfies the following properties:
\begin{enumerate}
\item $f(H) = 0$ if $H$ is trivial.
\item For every $H$ and every non-isolated vertex $v$ of $H$, $f(H\setminus N_H[v]) +1 \leqslant f(H)$.
\end{enumerate}
Then, for every subcomplex $\Gamma$ of $\Delta(G)$ with $\F(\Gamma) \subseteq \F(\Delta(G))$ we have $$\reg(I_\Gamma) \leqslant f(G)+1.$$
\end{cor}
\begin{proof} First we note that, for every graph $H$ and every $S\in\Delta(H)$ we have
$$\Delta(H\setminus N_H[S]) = \lk_{\Delta(H)}(S).$$
It implies that
$$\simp(\Delta(G)) = \{\Delta(H) \mid H\in \mathcal I_G\}.$$ 
Therefore, we can define a function $g \colon \simp(\Delta(G)) \to \N$, by sending $\Delta(H)$ to $f(H)$ for all $H\in \mathcal{I}_G$.

Note that for every graph $H$, we have $\Delta(H)$ is a simplex if and only if $H$ is trivial; and $\Delta(H)$ is a cone over a vertex $v$ if and only if $v$ is an isolated vertex of $H$. Together with the definition of the function $g$, it shows that $g$ satisfies all conditions of Theorem \ref{thm-delta}, and therefore by this theorem we obtain $\reg(I_\Gamma) \leqslant g(\Delta(G))+1 = f(G)+1$, as required.
\end{proof}

The theorem \ref{main-thm} when applying to an edge ideal of a graph has the following form.

\begin{lem}\label{lem-edge} Let $G$ be a graph. Then, 
$$\reg(I(G)^{(n)}) \leqslant 2(n-1) + b, \text{ for all } n\geqslant 1,$$ 
where $b = \max\{\reg(I_\Gamma)\mid \Gamma \text{ is a subcomplex of } \Delta(G) \text{ with } \F(\Gamma) \subseteq \F(\Delta(G))\}$.
\end{lem}
\begin{proof} Since $I(G) = I_{\Delta(G)}$ and $\delta(I(G)) = 2$ by \cite[Example 4.4]{DHHT}, therefore the lemma follows from Theorem \ref{main-thm}. 
\end{proof}
\medskip

We are now in position to prove the main result of this section.

\begin{thm}\label{EdgeBoundLinear} Let $G$ be a graph. Then,
$$\reg(I(G)^{(n)}) \leqslant 2n + \ordmatch(G)-1, \text{ for all } n\geqslant 1.$$
\end{thm}
\begin{proof} By Lemma \ref{lem-edge}, it remains to show that $\reg(I_\Gamma) \leqslant \ordmatch(G)+1$, for every subcomplex $\Gamma$ of $\Delta(G)$ with $\F(\Gamma) \subseteq \F(\Delta(G))$.

Consider the function $f\colon \mathcal{I}_G \to \N$ defined by 
$$f(H)=
\begin{cases}
0 & \text{ if } H \text{ is trivial},\\
\ordmatch(H) & \text{ otherwise}.
\end{cases}
$$
For every non-isolated vertex $v$ of $H$, we have $f(H\setminus N_H[v]) +1 \leqslant f(H)$ by \cite[Lemma 2.1]{Fv}, hence $f$ satisfies all conditions of Corollary \ref{lem-f}, so that by this corollary
$$\reg(I_\Gamma) \leqslant f(G)+1 = \ordmatch(G)+1,$$
and the theorem follows.
\end{proof}

\begin{rem} Let $G$ be a graph with $\ordmatch(G) = \nu(G)$. Then,
$$\reg(I(G)^{(n)}) = 2n + \nu(G)-1, \text{ for all } n\geqslant 1.$$

Indeed, for every positive integer $n$, the lower bound $\reg(I(G)^{(n)}) \geqslant 2n + \nu(G)-1$ comes from Lemma \ref{ind-reg}, and the upper bound follows from Theorem \ref{EdgeBoundLinear} because $\ordmatch(G) = \nu(G)$.

As a consequence, we quickly recover the main result of Fakhari in \cite{Fc}, which says that the equality holds when $G$ is a {\it Cameron-Walker} graph, where a graph $G$ is called Cameron-Walker if $\nu(G) = \match(G)$ (see e.g. \cite{HHKO}).  For such a graph $G$,  $\ordmatch(G) = \nu(G)$ since $\nu(G) \leqslant \ordmatch(G) \leqslant  \match(G)$.\end{rem}

\subsection*{Acknowledgment} We are supported by Project ICRTM.02\_2021.02 of the International Centre for Research and Postgraduate Training in Mathematics (ICRTM), Institute of Mathematics, VAST.

\end{document}